\documentclass{amsart}

\usepackage{amssymb}
\usepackage{amscd}
\usepackage{latexsym}
\usepackage{amsfonts}
\usepackage{amsmath}

\newtheorem{proposition}{Proposition}%[section]
\newtheorem{lemma}{Lemma}%[section]
\newtheorem{theorem}{Theorem}%[section]
\newtheorem{definition}{Definition}%[section]
\newcommand{\eproof}{\begin{flushright} $\square$ \end{flushright}}

\newcommand{\im}{\mathop{\fam0 Im}\nolimits}

\newcommand{\tr}{\mathop{\fam0 Tr}\nolimits}
\newcommand{\re}{\mathop{\fam0 Re}\nolimits}

\newcommand{\Hom}{\mathop{\fam0 Hom}\nolimits}
\newcommand{\End}{\mathop{\fam0 End}\nolimits}

\newcommand{\ad}{\mathop{\fam0 ad}\nolimits}

\newcommand{\Id}{\mathop{\fam0 Id}\nolimits}

\newcommand{\Pic}{\mathop{\fam0 Pic}\nolimits}

\newcommand{\bC}{\mathop{\fam0 {\mathbb C}}\nolimits}

\newcommand{\C}{C}

\newcommand{\Z}{\mathop{\fam0 {\mathbb Z}}\nolimits}
\newcommand{\N}{\mathop{\fam0 {\mathbb N}}\nolimits}
\newcommand{\bZ}{\Z{}}

\newcommand{\ra}{\mathop{\fam0 \rightarrow}\nolimits}

\newcommand{\cD}{\mathop{\fam0 {\mathcal D}}\nolimits}
\newcommand{\cV}{\mathop{\fam0 {\tilde {\mathcal V}}}\nolimits}
\newcommand{\cH}{\mathop{\fam0 {\mathcal H}}\nolimits}
\newcommand{\tcH}{\mathop{\fam0 {\tilde {\mathcal H}}}\nolimits}
\newcommand{\T}{\mathop{\fam0 {\mathcal T}}\nolimits}
\newcommand{\E}{\mathop{\fam0 {\mathcal E}}\nolimits}

\newcommand{\fg}{\mathop{\fam0 {\mathfrak g}}\nolimits}
\newcommand{\V}{\mathop{\fam0 {\mathcal V}}\nolimits}
\newcommand{\s}{\sigma}
\renewcommand{\L}{\mathop{\fam0 {\mathcal L}}\nolimits}

\renewcommand{\P}{\mathop{\fam0 {\mathbb P}}\nolimits}

\newcommand{\Nabla}{\mathop{\fam0 {\mathbf {\hat \nabla}}}\nolimits}
\newcommand{\Nablae}{\mathop{\fam0 {\mathbf {\hat \nabla}}^e}\nolimits}

\newcommand{\tOmega}{\mathop{\fam0 {{{\tilde \Omega}}}}\nolimits}
\newcommand{\tPhi}{\mathop{\fam0 {{{\tilde \Phi}}}}\nolimits}
\newcommand{\bH}{\mathop{\fam0 {{\mathbb H}}}\nolimits}

\newcommand{\Teim}{Teichm{\"u}ller }
\newcommand{\Ric}{\mathop{\fam0 Ric}\nolimits}

\begin{document}

\title[Mapping Class Groups do not have Kazhdan's Property (T)]
{Mapping Class Groups do not have\\ Kazhdan's Property (T)}

%\date{15'th of April 2002}
\author{J\o rgen Ellegaard Andersen}
\address{Department of Mathematics\\
        University of Aarhus\\
        DK-8000, Denmark}
\email{andersen@imf.au.dk}

\begin{abstract}
We prove that the mapping class group of a closed oriented surface
of genus at least two does not have Kazhdan's property (T).

\end{abstract}

\maketitle

\section{Introduction}

In the paper \cite{Kazh} Kazhdan introduced his property (T) for topological groups.
A topological group has Kazhdan's property (T), if the trivial representation is
isolated in the Fell topology on the space of unitary representations of the group.
Alternatively, we can formula Kazhdan's property (T) as follows.

\begin{definition}[Kazhdan]
For a unitary Hilbert space representation $\rho$ of a topological group $G$, a unit vector $v$
is call $(\epsilon,K)$-invariant, where $\epsilon$ is positive real and $K$ is a compact
subset of $G$, if
$$ \mid \rho(g)v -v \mid <\epsilon \mbox{, } \forall g\in K.$$
We say that $\rho$ has almost invariant vectors if there exist $(\epsilon,K)$-invariant
unit vector for all such pairs $(\epsilon,K)$.
A topological group $G$ has Kazhdan's property (T) if every unitary Hilbert space representation
which has almost invariant vectors has an actual nontrivial invariant vector.
\end{definition}

Since Kazhdan introduced this property it has been rather extensively studied, also
for discrete countable groups, as we shall be interested here. For such groups we
can give an alternative formulation of Kazhdan's property (T).

\begin{definition}
Let $\rho$ be a unitary Hilbert space representation of a discrete countable group $G$.
By an almost fixed vector for $\rho$ we mean a sequence of unit vectors
$(v_k)\subset \cH$ with the property that
$$\lim_{k \ra \infty} \mid \rho(g)v_k -v_k \mid =0$$
for all $g\in G$.
\end{definition}

We see that a discrete countable group $G$ has Kazhdan's property (T) if and only if
the existence of an
almost fixed vector implies the existence of a non trivial fixed vector for
all unitary Hilbert space representations of $G$.

\begin{theorem}\label{Main}
The mapping class group of a closed oriented surface of genus at
least two does not have Kazhdan's property (T).
\end{theorem}

We construct a counter example to Kazhdan's property (T) for these
mapping class groups using the Reshetikhin-Turaev Topological
Quantum Field Theory constructed in \cite{RT1}, \cite{RT2} and
\cite{Tu}. These TQFT-constructions by Reshetikhin and Turaev was
given on the basis of the suggestions by Witten in \cite{W}, which
gave a quantum field theory description of the Jones polynomial
\cite{J}.

Indeed we shall need the geometric constructions of these TQFT's
as proposed by Witten in that paper and by Atiyah in \cite{At}.
That the geometric construction gives the same representations as
the Reshetikhin-Turaev TQFT representations follow from combining
the results of \cite{L} and \cite{AU1}, \cite{AU2}, \cite{AU3} and
\cite{AU4}. In fact this identifies the geometrically constructed
representations with the TQFT representations constructed by
Blanchet, Habegger, Masbaum and Vogel in \cite{BHMV1} and
\cite{BHMV2}, which is the skein theory construction of the
RT-TQFT's. Please see Theorem \ref{AUMain} below regarding this.

Let us briefly recall the geometric construction of these
representations of the mapping class group.

Let $\Sigma$ be a closed oriented surface of genus at least two.
Let $p$ be a point on $\Sigma$. Let $M$ be the moduli space of
flat $SU(2)$ connections on $\Sigma - p$ with holonomy around $p$
equal $-\Id \in SU(2)$. This moduli space is smooth and has a
natural symplectic structure $\omega$. There is natural smooth
symplectic action of the mapping class group $\Gamma$ of $\Sigma$
on $M$. More over there is a unique
prequantum line bundle $(\L,\nabla, (\cdot,\cdot))$ over
$(M,\omega)$. The Teichm\"{u}ller space $\T$ of complex structures on
$\Sigma$ naturally $\Gamma$-equivariantly parametrizes K\"{a}hler structures on $(M,\omega)$.
For $\sigma\in \T$, we denote $(M,\omega)$ with is corresponding
K\"{a}hler structure $M_\sigma$.

By applying geometric quantization to the moduli space $M$, one gets
a certain finite rank bundle over \Teim space $\T$, which we will
call the {\em Verlinde} bundle $\V_k$ at level $k$, where $k$ is
 any positive integer. The fiber of this
bundle over a point $\sigma\in \T$ is $\V_{k,\sigma} =
H^0(M_\sigma,\L^k)$. We observe that there is a natural Hermitian
structure $\langle\cdot,\cdot\rangle$ on
$H^0(M_\sigma,\L^k)$ by restricting the $L_2$-inner product
on global $L_2$ sections of $\L^k$ to $H^0(M_\sigma,\L^k)$.

The main result pertaining to this bundle $\V_k$ is that its
projectivization $\V_k$ supports a natural flat $\Gamma$-invariant connection
$\Nabla$. This is a result proved independently by Axelrod, Della
Pietra and Witten \cite{ADW} and by Hitchin \cite{H} (see also \cite{A3}). This flat
connection $\Nabla$ induces a {\em flat} connection $\Nablae$ in
$\End(\V_k)$. Let $\End_0(\V_k)$ be the subbundle consisting of
traceless endomorphisms. The connection $\Nablae$ also induces a
connection in $\End_0(\V_k)$, which is invariant under the action of $\Gamma$.

We get this way for each $k$ a finite dimensional representation
of $\Gamma$, namely on the covariant constant sections, say
$\cH_k$, of $\End_0(\V_k)$ over $\T$. Let
 $$\tcH = \bigoplus_{k+2 \text{ prime}}^\infty \cH_k$$
on which $\Gamma$ acts. From the proof of the asymptotic
faithfulness in \cite{A2}, one see that this representation of
$\Gamma$ is faithful.

Each of the vector spaces $\cH_k$ has a natural positive definite Hermitian structure
 $[ \cdot , \cdot ]$,
which is preserved by the action of $\Gamma$.
 This Hermitian structure is clear from the skein theory construction of $\cH_k$ following \cite{BHMV2}:

Consider the BHMV-TQFT (as defined in \cite{BHMV2}) at $A=
\exp(2\pi i/4(k+2))$. The label set for this theory is then $L_k =
\{0, 1, \ldots, k\}$. We denote by $Z_k$ the vector space this
theory associates to $\Sigma \sqcup \overline{\Sigma}$ with
$p\in\Sigma$ label by the last color $k$ in both copies of
$\Sigma$. The BHMV construction also requires us to choose a
$p_1$-structure on $\Sigma \sqcup \overline{\Sigma}$. However, we
note that this vector space does not depend on the choice of a
$p_1$-structure, as long as the choice of the $p_1$-structure
$\Sigma \sqcup \overline{\Sigma} = \partial(\Sigma \times [0,1])$
extends over $\Sigma \times [0,1]$.

Since the vector space $Z_k$ is part of a TQFT, there is in
particular an action of the mapping class group of $\Sigma\sqcup
\overline{\Sigma}$ on $Z_k$. There is a natural homomorphism of
$\Gamma$ into the mapping class group of $\Sigma\sqcup
\overline{\Sigma}$ given by mapping $\phi\in \Gamma$ to $\phi \sqcup
\phi$.

In \cite{BHMV2} a Hermitian structure $\{ \cdot , \cdot \}$ is
constructed on $Z_k$, which is invariant under the action of the
mapping class group of $\Sigma\sqcup
\overline{\Sigma}$ and therefore also invariant under the action
of $\Gamma$. For the choice of $A$ made above, it is proved in
\cite{BHMV2}, that the Hermitian structure $\{ \cdot , \cdot \}$
is positive definite.

By the work of Andersen and Ueno \cite{AU1}, \cite{AU2},
\cite{AU3} and \cite{AU4} combined with the work of Laszlo
\cite{L}, we have an identification of the two constructions.

\begin{theorem}[Andersen \& Ueno]\label{AUMain}
There is a natural $\Gamma$-equivariant isomorphism
\[I_k : Z_k \ra \cH_k\]
\end{theorem}

Using the isomorphism $I_k$, we define the the positive definite
Hermitian structure $[ \cdot , \cdot ]$ on $\cH_k$ by the
formula
\[[ \cdot , \cdot ] = \{ I^{-1}_k(\cdot) , I^{-1}_k(\cdot) \}.\]
The norm associated to $[ \cdot , \cdot ]$ is denoted $[ \cdot]$.
The Hermitian structure $[ \cdot , \cdot ]$ on $\cH_k$ induces
a Hermitian structure on $\End(\V_k)$, which is parallel with
respect to $\Nablae$ and which is $\Gamma$-invariant.

\begin{definition}
We define $\cH$ to be the Hilbert
space completion of $\tcH$ with respect to the norm  $[ \cdot ]$.
\end{definition}
This is an infinite
dimensional Hilbert space, on which $\Gamma$ acts isometrically.
This representation provides us with the needed counter example to
Kashdan's property (T) for $\Gamma$. Let us discuss the proof of
this.

\begin{theorem}[Roberts]\label{RobertsT}
The only $\Gamma$ invariant vector in $\cH$ is $0$.
\end{theorem}

This theorem follows from the fact that the representations
$\cH_k$ are irreducible as $\Gamma$-representation for $k$, such
that $k+2$ is prime. This results was established in the
un-twisted case by Roberts in \cite{Ro} and his proof can be
applied word for word also to this case.

The basic idea behind building the required almost fixed vector for
$\cH$ is to consider coherent states on $M_\sigma$, $\sigma \in \T$.

Fix a point $x\in M$. Evaluation at $x$ gives a section of $\V_k^*$
up to scale. Using $\langle\cdot,\cdot\rangle$ we get induced
a section $e_x^{(k)}$ of $\V_k$ up to scale. For each
$\sigma\in\T$,
$e_x^{(k)}(\sigma)$ is the coherent state associated to $x$ on $M_\sigma$.
Let $E_x^{(k)}$ be the section of
$\End(\V_k)$ obtained as the orthogonal projection (with respect to
$\langle\cdot,\cdot\rangle$) onto the one dimensional subspace
spanned by $e_x^{(k)}$. We observe that $E_x^{(k)}$ only
depends on $x$.

\begin{theorem}\label{ccE}
The sections $E_x^{(k)}$ of $\End(\V_k)$ over $\T$ are
asymptotically covariant constant. I.e. for any pair of points
$\sigma_0,\sigma_1 \in \T$ there exists a constant $C$ such that
$$\left| P^e_{\sigma_0,\sigma_1} \left( E_x^{(k)}(\sigma_0)\right)
- E_x^{(k)}(\sigma_1)\right| \leq \frac Ck,$$ where
$P^e_{\sigma_0,\sigma_1}$ is the parallel transport from
$\sigma_0$ to $\sigma_1$ in $\End(\V_k)$ and the norm $|\cdot|$ is
the one associated to the Hermitian structure on $\V_k^* \otimes
\V_k$ induced from $\langle\cdot,\cdot\rangle$ on $\V_k$.
\end{theorem}
The proof of this theorem is given in section \ref{cscc}.

In order to relate these norm estimates to estimates for the norm
associated to the Hermitian structure $[\cdot,\cdot]$, we need the
following result.

Pick a point $\sigma_0 \in \T$.

\begin{theorem}[Andersen]\label{normeq}
There exist a constant $C$ such that at $\sigma_0$
\[ C^{-1} | \Psi| \leq [\Psi] \leq C  | \Psi|\]
for all $\Psi$ in $\End(\V_k)_{\sigma_0}$ and all $k\in \N$.
\end{theorem}

This theorem is proved in \cite{A6}.

To produce the almost fixed vectors, we now pick a finite
subgroup $\Lambda$ of $SU(2)$ which contains $-1\in SU(2)$ and we consider the
finite subset $X$ of $M$, consisting of connections which reduces to $\Lambda$.
As we will see in section \ref{PF} for appropriate choice of $\Lambda$, we
get a non empty finite subset of $M$ this way which is
invariant under the action of the mapping class group and such that $|X|>1$.
Let now $E_X^{(k)}$ be the
section of $\End(\V_k)$ given by
$$E_X^{(k)} = \sum_{x\in X} E_x^{(k)}.$$
Let $E_{X,0}^{(k)}$ be the
traceless part of $E_X^{(k)}$. As we will see in section \ref{PF}, for large enough $k$,
$E_X^{(k)} \neq 0$. Hence for large enough $k$ we have a unique
vector in $\E_{X,0}^{(k)}\in \cH_k$, which at $\sigma_0$ agrees with
$E_{X,0}^{(k)}(\sigma_0)/[E_{X,0}^{(k)}(\sigma_0)]$.

\begin{theorem}\label{AFV}
The sequence $\{\E_{X,0}^{(k)}\}$
is an almost fixed vector for the action of $\Gamma$ on $\cH$.
\end{theorem}

This theorem will be proved in section \ref{PF}. Our main Theorem
\ref{Main}
is of course a consequence of Theorem \ref{RobertsT} and Theorem
\ref{AFV}.

Since the mapping class group of a genus one closed surface is
$SL(2,\Z)$, it is well known that this mapping class group does
not have Kazhdan's property (T). The constructions as presented
here can be applied to the $U(1)$-moduli space in the genus one
case to provide a counter example for $SL(2,\Z)$ as well.

It is a result of F. Taherkhani
that the mapping class group of a closed oriented surface of genus two does not have
Kazhdan's property (T) \cite{Ta}. Taherkhani's proof of the genus two case relies on computer
aided calculations of the first cohomology of certain cofinite subgroups of the genus two mapping group
and so is completely different from the arguments presented here.

The construction of the almost fixed vector applies to the rather more general
setting of \cite{A3}, for which we have established the existence
of Hitchin's connection. It would be interesting to find a
geometric argument for the irreducibility, which potentially could be
applied in the more general setting of \cite{A3}, and thereby possibly provide
counter examples to Kazhdan's property (T) for other groups.

The methods of asymptotic analysis and the theory of Toeplitz operators
has also allowed us to link these TQFT to the Nielsen-Thurston classification
of mapping classes \cite{A4} and to asymptotics of Hermitian pairings of
loop operators \cite{A5} (Toeplitz operator interpretation of \cite{MN}).

The paper is organized as follows. In section \ref{sec2} and
\ref{sec3} we recall the construction of the Verlinde bundle and
Hitchin's connection. Basics about coherent states and Toeplitz
operators are recalled in section \ref{sec4}. Another
geometrically defined Hermitian structure which is asymptotically
preserved by the Hitchin connection is recalled in section
\ref{Hermstruc}. We need this structure in section \ref{cscc} to
prove the asymptotic flatness of the coherent states as stated in
Theorem \ref{ccE}, followed in section \ref{PF} by the
construction of the almost fixed vector.

We would like to thank Vaughan Jones for suggesting the
use of these representations to settle this property for these
mapping class groups. Jones has also posted this problem on the
CTQM problem list at www.ctqm.au.dk (Problem 14). At the time
of this writing it remains an open problem
to what extend the mapping class groups has the Haagerup property. Further we would like to thank Gregor Masbaum
for valuable discussions.

\section{The gauge theory construction of the Verlinde bundle}\label{sec2}

Let us now very briefly recall the construction of the Verlinde
bundle. Only the details needed in this paper will be given. We
refer e.g. to \cite{H} for further details. As in the introduction
we let $\Sigma$ be a closed oriented surface of genus $g\geq 2$
and $p\in \Sigma$. Let $P$ be a principal $SU(2)$-bundle over
$\Sigma$. Clearly, all such $P$ are trivializable. Let $M$ be the
moduli space of flat $SU(2)$-connections in $P|_{\Sigma - p}$ with
holonomy $-1 \in SU(2)$ around $p$. We can identify
\[M = \Hom'({\tilde \pi}_1(\Sigma), SU(2) )/ SU(2).\]
Here ${\tilde \pi}_1(\Sigma)$ is the universal central extension
\[0 \ra \bZ \ra {\tilde \pi}_1(\Sigma) \ra \pi_1(\Sigma) \ra 1\]
as discussed in \cite{H} and in \cite{AB} and $\Hom'$ means the
space of homomorphisms from ${\tilde \pi}_1(\Sigma)$ to $SU(2)$
which send the image of $1\in \bZ$ in ${\tilde \pi}_1(\Sigma)$ to
$-1\in SU(2)$ (see \cite{H}).

\begin{theorem}[Goldman; Atiyah \& Bott]
The moduli space $M$ is a smooth compact manifold of dimension
$6g-6$.
\end{theorem}

See \cite{G} for a representation variety proof of this theorem
and \cite{AB} for a gauge theory proof.

Since there is a natural homomorphism from the mapping class group
to the outer automorphisms of ${\tilde \pi}_1(\Sigma)$, we get a
smooth action of $\Gamma$ on $M$.

On $\fg = \mbox{Lie}(SU(2))$ we have the invariant symmetric
bilinear form $(X,Y) \mapsto \tr(XY)$, normalized such that
$-\frac{1}{6}\tr(\vartheta\wedge [\vartheta\wedge\vartheta])$ is a
generator of the image of the integer cohomology in the real
cohomology in degree $3$ of $SU(2)$, where $\vartheta$ is the
$\fg$-valued Maurer-Cartan $1$-form on $SU(2)$. This bilinear form
induces a symplectic form on $M$. At a flat connection $A$
representing a point $[A]\in M$:
\[\omega(\varphi_1,\varphi_2) = \int_{\Sigma}\tr(\varphi_1\wedge\varphi_2),\]
where $\varphi_i$ are $d_A$-closed $1$-forms on $\Sigma$ with
values in $\ad P$ representing tangent vectors to $M$ at $[A]$.
See e.g. \cite{G}, \cite{AB} or \cite{H} for further details on
this. We summarize this in the following theorem.

\begin{theorem}[Goldman; Atiyah \& Bott; Narashimhan \& Seshadri]
On the moduli space $M$, the form $\omega$ is a symplectic
structure and the natural action of $\Gamma$ on $M$ is symplectic.
\end{theorem}

Let $\L$ be the Hermitian line bundle over $M$ and $\nabla$ the
compatible connection in $\L$ constructed by Freed in \cite{Fr}.
This is the content of Corollary 5.22, Proposition 5.24 and
equation (5.26) in \cite{Fr} (see also the work of Ramadas, Singer
and Weitsman
\cite{RSW}). By Proposition 5.27 in \cite{Fr} we have that the
curvature of $\nabla$ is $\frac{\sqrt{-1}}{2\pi}\omega$. We will also use
the notation $\nabla$ for the induced connection in $\L^k$, where
$k$ is any integer.

\begin{theorem}[Ramadas, Singer \& Weitsman; Freed]
The Hermitian line bundle with connection $(\L,\nabla)$ is a
prequantum line bundle over the moduli space, i.e. the curvature
of $\nabla$ is $\frac{\sqrt{-1}}{2\pi}\omega$.
\end{theorem}

By an almost identical construction, we can lift the action of
$\Gamma$ on $M$ to act on $\L$ such that the Hermitian connection
is preserved (See e.g. \cite{A1}). In fact, since $H^2(M,
\bZ)\cong \bZ$ and $H^1(M,\bZ) = 0$, it is clear that the action
of $\Gamma$ leaves the isomorphism class of $(\L,\nabla)$
invariant, thus alone from this one can conclude that a $U(1)$-central
extension of $\Gamma$ acts on $(\L,\nabla)$ covering the $\Gamma$
action on $M$. This is actually all we need in this paper. We will
return to this point at the end of this section.

Let now $\sigma\in \T$ be a complex structure on $\Sigma$. Let us
review how $\sigma$ induces a complex structure on $M$ which is
compatible with the symplectic structure on this moduli space. The
complex structure $\s$ induces a $*$-operator on $1$-forms on
$\Sigma$ with values in $\ad P$, which acts on the harmonic forms
with square $-1$. Hence we get an almost complex structure on $M$
by letting $I= I_\sigma = -*$ acting on harmonic forms $\Sigma$
with values in $\ad P$.

We have the following classical result by Narasimhan and Seshadri
(see \cite{NS}),

\begin{theorem}[Narasimhan \& Seshadri]
The triple $(M,\omega, I_\sigma)$ is a smooth K{\"a}hler manifold
for any $\s\in \T$.
\end{theorem}

We use the notation $M_\sigma = (M,\omega, I_\sigma)$. By using
the $(0,1)$ part of $\nabla$ in $\L$ over $M_\s$, we get an
induced holomorphic structure in the bundle $\L$. See also \cite{H}
and \cite{AB} for further details on this.

From a more algebraic geometric point of view, we consider the
moduli space of S-equivalence classes of semi-stable bundles of
rank $2$ and determinant isomorphic to the line bundle ${\mathcal
O}([p])$. By using Mumford's Geometric Invariant Theory,
Narasimhan and Seshadri (see \cite{NS}) showed that

\begin{theorem}[Narasimhan \& Seshadri]
The moduli space moduli space of semi-stable bundles of rank $2$
and determinant isomorphic to ${\mathcal O}([p])$ is a smooth
complex algebraic projective variety, which is isomorphic as a
K{\"a}hler manifold to $M_\sigma$.
\end{theorem}

Referring to \cite{DN} we further recall that

\begin{theorem}[Drezet \& Narasimhan]\label{DNTh1}
The Picard group of $M_\sigma$ is generated by the holomorphic line bundle
$\L$ over $M_\sigma$
constructed above:
\[\Pic(M_\sigma) = \langle \L \rangle. \]
\end{theorem}

\begin{definition}\label{Verlinde}
The {\em Verlinde bundle} $\V_k$ over \Teim space is by definition the
bundle whose fiber over $\sigma\in \T$ is $H^0(M_\sigma,\L^k)$, where $k$ is a
positive integer.
\end{definition}

We will consider the endomorphism bundle $\End(\V_k)$ of $\V_k$.
We observe that the general argument above gives an action of a
central extension of $\Gamma$ acting on $\V_k$, which then induces
an action of $\Gamma$ on $\End(\V_k)$. This intern induces an
action of $\Gamma$ on the sub-bundle $\End_0(\V_k)$ consisting of
traceless endomorphisms.

\section{The projectively flat connection}\label{sec3}
In this section we will review Axelrod, Della Pietra and Witten's
and Hitchin's construction of the projective flat connection over
\Teim space in the Verlinde bundle.
We refer to \cite{H},
\cite{ADW} and \cite{A3}
for further details.

Let $\cV_k$ be the trivial $C^\infty(M,\L^k)$-bundle over $\T$
which contains $\V_k$, the Verlinde sub-bundle.

The Hitchin connection is a connection in $\cV_k$ which is of the form
\begin{equation}\Nabla_{v} = \Nabla^t_{v} - u(v),\label{Con}
\end{equation}
for all $v\in T(\T)$. Here $\Nabla^t$ is the trivial connection
in the trivial bundle $\cV_k$ and $u$ is a smooth map from $T(\T)$
to the vector space $D(M,\L^k)$ consisting of differential
operators acting on sections of $\L^k$. Hitchin constructs a
specific $u$ such that the corresponding connection preserved the
Verlinde subbundle. Let us recall his construction here.

The holomorphic tangent space to \Teim space at $\sigma\in \T$ is
given by
\[T^{1,0}_\sigma (\T) \cong H^1(\Sigma_\sigma, K_\sigma^{-1}).\]
If $v\in T_\sigma(\T)$, then we denote its $(1,0)$-part in $T^{1,0}_\sigma (\T)$
by $v'$. The holomorphic co-tangent space to the moduli space of
semi-stable bundles at the equivalence class of a stable bundle
$E$ is given by
\[T^*_{[E]}M_\sigma \cong H^0(\Sigma_\sigma, \End_0(E)\otimes K_\sigma).\]
We denote the holomorphic tangent bundle of $M_\sigma$ by $T_\sigma$.
For each $v\in T_\sigma(\T)$ we now specify a
$G(v)\in \Omega^0(M_\sigma,S^2(T_\sigma))$ as a quadratic function on
$T^*_\sigma$ by the formula
\[G(v)(\alpha,\alpha) = \int_{\Sigma}\tr(\alpha^2)v',\]
for all $\alpha \in H^0(\Sigma, \End_0(E)\otimes K_\sigma)$.
From this formula it is
clear that $G(v)\in H^0(M_\s,S^2(T_\sigma))$.

Axelrod, Della Pietra and Witten's $u(v)$, which by the results of \cite{A3}
 agrees projectively with Hitchin's $u(v)$, given in
\cite{H}, is
\begin{equation}
u(v)(s) = \frac{1}{2(k + 2)}(\Delta_{G(v)} - 2
\nabla_{G(v) \partial F(\sigma)} + k v'[F](\sigma)) s.
\label{Hitchcon}
\end{equation}

The leading order term $\Delta_{G(v)}$ is the 2'nd order operator given by
\[
\begin{CD}
\Delta_{G(v)} : \C^\infty(M,\L^k) @>{\nabla^{1,0}}>>
\C^\infty(M,T_\sigma^*\otimes\L^k)
@>{G(v)}>> \C^\infty(M,T_\sigma\otimes\L^k)\\
 @>{\nabla^{1,0}\otimes 1 + 1\otimes\nabla^{1,0}}>>
 \C^\infty(M,T_\sigma^*\otimes T_\sigma\otimes \L^k) @>{\tr}>>
 \C^\infty(M,\L^k),
\end{CD}
\]
where we have used the Chern connection in $T_\sigma$ on the K\"{a}hler
manifold $(M_\sigma,\omega)$.

Further $F$ is the unique smooth map from $\T$ to the vector space
$C^\infty(M)$, such that $F(\sigma)$ is the Ricci potential
uniquely determined as the real
function with zero average over $M$, which satisfies the following
equation
\begin{equation}
\Ric_\sigma = 2 n \omega + 2 \sqrt{-1} \partial\bar{\partial} F(\sigma).
\label{riccipot}
\end{equation}

The notation $v'[F](\sigma)\in C^\infty(M)$ mean the derivative of $F$ along the
direction $v'$ in $\T$ at $\sigma\in \T$.

The complex vector field $G(v)\partial F(\sigma) \in \C^\infty(M, T_\sigma)$ is simply just the
contraction of $G(v)$ with $\partial F(\sigma) \in \C^\infty(M, T_\sigma^*)$.

We observe that
$\Nabla$ agrees with $\Nabla^t$ along the anti-holomorphic
directions $T^{0,1}(\T)$.

We remark that
there is some finite set of vector
fields  $X_r(v), Y_r(v), Z(v) \in
\C^\infty(M_\sigma,T)$, $r = 1, \ldots, R$ (where $v\in T_\sigma(\T)$) all varying
smoothly\footnote{This makes sense when we consider the
holomorphic tangent bundle $T$ of $M_\sigma$ inside the
complexified real tangent
bundle $TM \otimes {\mathbb C}$ of $M$.} with
$v\in T(\T)$, such that
\begin{equation}
\Delta_{G(v)} -2 \nabla_{G(v)\partial F} = \sum_{r=1}^R \nabla_{X_r(v)}
\nabla_{Y_r(v)} + \nabla_{Z(v)}.\label{Opform}
\end{equation}
This follows immediately from the definition of $\Delta_{G(v)}$.
See also \cite{A3}.

\begin{theorem}[Axelrod, Della Pietra \& Witten; Hitchin]\label{Pflat}
The expression (\ref{Con}) above defines a connection $\Nabla$ in the bundle
$\V_k$, which induces a flat connection in $\P(\V_k)$.
\end{theorem}

We remark about genus $2$, that \cite{ADW} covers this case, but
\cite{H} excludes this case, however, the work of Van Geemen and De
Jong \cite{vGdJ} extends Hitchin's approach to the genus $2$ case.
See also \cite{A3} where the connection is constructed in
a more general situation covering this moduli space application for all
$g>1$.
It is clear from formula (\ref{Con}) that $\Nabla$ is invariant
under the action of $\Gamma$ on $\P(\V_k)$.

We  are here interested in the induced {\em flat}
connection $\Nablae$ in the endomorphism bundle $\End(\V_k)$.
Suppose $\Phi$ is a section of $\End(\V_k)$. Then for all sections
$s$ of $\V_k$ and all $v\in T(\T)$ we have that
\[(\Nablae_v \Phi) (s) = \Nabla_v \Phi(s) - \Phi(\Nabla_v(s)).\]

It is clear from the construction of $\Nablae$, that the subbundle
$\End_0(\V_k)$ is preserved by $\Nablae$. Thus
$(\End_0(\V_k),\Nablae)$ is a vector bundle over $\T$ with a flat
mapping class group invariant connection.

\begin{definition}
Let $\cH_k$ be the finite dimensional representation of $\Gamma$ consisting
of the
covariant constant sections  of $(\End_0(\V_k),\Nablae)$ over
$\T$.
\end{definition}

The dimension of $\cH_k$ is of course $d_g(k)^2 - 1$, where
$d_g(k)$ is the rank of $\V_k$, which is given by the twisted
Verlinde formula \cite{Th} (see also \cite{AM})
\[d_g(k) = (k+1)^{g-1} \sum_{j=1}^{2k+1} (-1)^{j+1}
\left( \sin(\frac{\pi j}{2(k+1)})\right)^{2-2g} \]

\section{Asymptotics of  Toeplitz operators and coherent states}\label{sec4}

We shall in this section discuss the asymptotics of Toeplitz
operators and coherent states as the level $k$ goes to infinity.
The properties we need can all be derived from the fundamental
work of Boutet de Monvel and Sj\"{o}strand. In \cite{BdMS} they did a
microlocal analysis of the Szeg\"{o} projection, which can be applied
to the asymptotic analysis in the situation at hand, as it was
done by Boutet de Monvel and Guillemin in \cite{BdMG} (in fact in
a much more general situation than the one we consider here) and
others following them. In particular the applications developed by
Schlichenmaier \cite{Sch}, \cite{Sch1}, \cite{Sch2} and further by
Karabegov and Schlichenmaier \cite{KS} to the study of Toeplitz
operators in the geometric quantization setting, is what will
interest us here. Let us first describe the basis setting.

On $\C^\infty(M,\L^k)$ we have the $L_2$-inner product:
\[\langle s_1, s_2 \rangle = \frac{1}{m!}\int_M (s_1,s_2) \omega^m\]
where $s_1, s_2 \in \C^\infty(M,\L^k)$ and $(\cdot,\cdot)$ is the fiberwise
Hermitian structure in $L^k$.

Let $\sigma$ be a point in $\T$ and consider the K\"{a}hler manifold $M_\sigma$.
Inside the space of all smooth sections
$\C^\infty(M,\L^k)$, we have the finite dimensional subspace $H^0(M_\sigma ,\L^k)$ consisting
of holomorphic sections with respect to the complex structure on $M_\sigma$.

The $L_2$-inner product on $\C^\infty(M,\L^k)$ determines the
 orthogonal projection $\pi^{(k)}_\sigma : \C^\infty(M,\L^k) \ra H^0(M_\sigma ,\L^k)$.
For each $f\in \C^\infty(M)$ consider the associated {\em Toeplitz
operator} $T_{f,\sigma}^{(k)}$ given as the composition of the
multiplication operator (which we also denote f) $f :
\C^\infty(M,\L^k) \ra \C^\infty(M,\L^k)$ with the orthogonal
projection $\pi^{(k)}_\sigma : \C^\infty(M,\L^k) \ra
H^0(M_\sigma,\L^k)$, so that
\[T_{f,\sigma}^{(k)}(s) = \pi^{(k)}_\sigma(f s).\]
Since the multiplication operator is a zero order differential
operator, $T_{f,\sigma}^{(k)}$ is a zero-order Toeplitz operator.

The first result we will need is
due to Bordemann, Meinrenken and Schlichenmaier (see \cite{BMS}).
The $L_2$-inner product on $\C^\infty(M,\L^k)$ induces an inner product $\langle\cdot,\cdot\rangle$ on $H^0(M_\sigma,\L^k)$, which in turn induces the operator norm $\|\cdot\|$ on $\End(H^0(M_\sigma,\L^k))$.

\begin{theorem}[Bordemann, Meinrenken and Schlichenmaier]\label{BMS1}
For any $f\in \C^\infty(M)$ we have that
\[\lim_{k\ra \infty}\|T_{f,\sigma}^{(k)}\| = \sup_{x\in M}|f(x)|.\]
\end{theorem}

This result follows also directly from the results of \cite{BdMS} and \cite{BdMG} as shown in
\cite{KS}.

\begin{theorem}[Schlichenmaier]\label{S}
For each $\sigma\in \T$ and any pair of smooth functions $f_1,
f_2\in \C^\infty(M)$, we have an asymptotic expansion
\[T_{f_1,\sigma}^{(k)}T_{f_2,\sigma}^{(k)} \sim \sum_{l=0}^\infty
T_{c^{(l)}_\sigma(f_1,f_2),\sigma}^{(k)} k^{-l},\] where
$c^{(l)}_\sigma(f_1,f_2) \in C^\infty(N)$ are uniquely determined
since $\sim$ means the following: For all $L\in \Z_+$ we have that
\begin{equation}
\|T_{f_1,\sigma}^{(k)}T_{f_2,\sigma}^{(k)} - \sum_{l=0}^L
T_{c^{(l)}_\sigma(f_1,f_2),\sigma}^{(k)} k^{-l}\| =
O(k^{-(L+1)}).\label{normasympToep}
\end{equation}
Moreover, $c^{(0)}_\sigma(f_1,f_2) = f_1f_2$.
\end{theorem}

This Theorem was proved in \cite{Sch} (again building on the works
\cite{BdMS} and \cite{BdMG}) and it is published in \cite{Sch1}
and \cite{Sch2}, where it is also proved that the formal
generating series for the $c_l(f_1,f_2)$'s gives a formal
deformation quantization of the Poisson structure on $(M,\omega)$.
By examining the proof in \cite{Sch} (or in \cite{Sch1} and
\cite{Sch2}) of this Theorem, one observes that for continuous
families of functions, the estimates in Theorem \ref{S} are
uniform over compact subsets of $\T$.

Let $\sigma_0$ and $\sigma_1$ be two points in $\T$. For any $f\in \C^\infty(M)$
we consider
$$T^{(k)}_{f,(\sigma_0,\sigma_1)} : H^0(M_{\sigma_0},\L^k) \ra H^0(M_{\sigma_1},\L^k)$$
given by
$$ T^{(k)}_{f,(\sigma_0,\sigma_1)}  = T^{(k)}_{f,\sigma_1} \mid_{H^0(M_{\sigma_0},\L^k)}.$$
We see that
$$ T^{(k)}_{f,(\sigma_0,\sigma_1)} = \pi^{(k)}_{\sigma_1} f \pi^{(k)}_{\sigma_0}.$$
We will also use the notation
$$\pi^{(k)}_{(\sigma_0,\sigma_1)} = \pi^{(k)}_{\sigma_1}\mid_{H^0(M_{\sigma_0},\L^k)}.$$

\begin{theorem}\label{TOA2}
For all $f\in \C^\infty(M)$ we have that
$$ \| T^{(k)}_{f,(\sigma_0,\sigma_1)} - \pi^{(k)}_{(\sigma_0,\sigma_1)} T^{(k)}_{f,\sigma_0}\| = O(k^{-1})$$
and
$$ \| T^{(k)}_{f,(\sigma_0,\sigma_1)} -  T^{(k)}_{f,\sigma_1}\pi^{(k)}_{(\sigma_0,\sigma_1)}\| = O(k^{-1}).$$
\end{theorem}

We will prove this theorem by using the theory of Fourier integral
operators and their symbol calculus as discussed in \cite{BdMG}.

Let $Z$ be the unit tangent bundle in $\L^*$. Let $\tOmega$ be the
volume form on $Z$, which is invariant under the $U(1)$ action on
$Z$ and with the property that
\[\int_Z \tau^*(f) \tOmega = \int_M f\Omega\]
for all $f\in C^\infty(M)$. The vector space $C^\infty(M,L^k)$ is
isommorphic to the vector space of smooth functions on $Z$ which
transforms in the $k$'th representation of $U(1)$ via the formula
\[\psi_s(\alpha) = \alpha^{\otimes k}(s(x))\]
for all $\alpha\in \tau^{-1}(x)$, $x\in M$ and all $s\in
C^\infty(M,L^k)$. In fact this is an isometry between $L_2(M,L^k)$
and the $k$'th weight space of the $U(1)$ action on $L_2(Z)$.

For each $\sigma\in \T$, we define the Hardy space $\bH_\sigma
\subset L_2(Z)$ consisting of the $L_2$-functions, which extend
over the unit disc bundle in $\L^*$ holomorphically with respect
to the complex structure induced from the one on $M_\sigma$. The
orthogonal project from $L_2(Z)$ to the closed subspace
$\bH_\sigma$ is the Szeg\"{o} projections and we denote it
$\Pi_\sigma$. If we denote by $\bH_\sigma^{(k)}$ the $k$'th weight
space of the $U(1)$ on $\bH_\sigma$. Then the above isomorphism
restrict to an isometry $H^0(M_\sigma, \L^k) \cong
\bH_\sigma^{(k)}$ and we (also) denote the Bargman orthogonal
projection onto $\bH_\sigma^{(k)}$ by $\pi_\sigma^{(k)}$.

This way, we see that every $\sigma\in \T$ makes $Z$ a
pseudo-convex domain with a Toeplitz structure $\Pi_\sigma$ in the
sense of \cite{BdMG}.

 \proof[Proof of Theorem \ref{TOA2}]
 The Szeg\"{o}
projectors $\Pi_{\sigma_0}$ and $\Pi_{\sigma_1}$ with respect to
the complex structure $\sigma_0$ and $\sigma_1$ are Fourier
integral operators by Theorem 11.1 in \cite{BdMG} and by the
composition rule of symbols of these types of operators (Theorem
9.8 in \cite{BdMG}), we see that the leading order symbol of the
operator
 $$\Pi_{\sigma_1} f \Pi_{\sigma_0} - \Pi_{\sigma_1} \Pi_{\sigma_0} f \Pi_{\sigma_0}$$
vanishes. This means this zero order operator in fact is of order
$-1$. Consider the first order differential operator $D$
corresponding to the infinitesimal generator of the circle action
on $X$. This operator commutes with $\Pi_{\sigma_1}$ and we have
that
 $$\Pi_{\sigma_1} Df \Pi_{\sigma_0} - \Pi_{\sigma_1} D\Pi_{\sigma_0} f \Pi_{\sigma_0}$$
is a zero order operator, hence bounded. However $D$ acts on
$H^0(M_{\sigma_1},\L^k)$ by multiplication by $k$. Thus there
exist a constant $C$ (equal to the operator norm of the above
operator) such that
 $$ \| k (\pi^{(k)}_{\sigma_1} f \pi^{(k)}_{\sigma_0} - \pi^{(k)}_{\sigma_1} \pi^{(k)}_{\sigma_0} f \pi^{(k)}_{\sigma_0}) \| \leq C.$$
The other inequality is proved the same way. \eproof

We again observe from the proof, that the estimates in Theorem
\ref{TOA2} is uniform for $(\sigma_0,\sigma_1)$ contained in
compact subsets of $\T$.

Let us now consider coherent states. Pick a point $x$ in $M$ and
let $\alpha \in \L^*_x$. We think of $\alpha$ as a linear map
\[\alpha : \C^\infty(M,\L^k) \ra \bC\]
given by
\[\alpha(s) = \alpha^{\otimes k}(s(x)).\]

For a $\sigma$ in $\T$, we let $e^{(k)}_{\alpha,\sigma}\in
H^0(M_{\sigma},\L^k)$ be the corresponding coherent state, i.e.
\begin{equation}
\langle s, e^{(k)}_{\alpha,\sigma} \rangle = \alpha(s)
\label{defcs}
\end{equation}
for all $s\in H^0(M_\sigma,\L^k)$.

\begin{theorem}\label{TOCSA}
For all $f\in \C^\infty(M)$, any $x\in M$ and any pair of points
$(\sigma_0,\sigma_1)$ in $\T$ we have that
 $$ \left|
T^{(k)}_{f,(\sigma_0,\sigma_1)} \frac{e^{(k)}_{\alpha,\sigma_0}
}{|e^{(k)}_{\alpha,\sigma_0}|}- f(x)
\frac{e^{(k)}_{\alpha,\sigma_1}
}{|e^{(k)}_{\alpha,\sigma_1}|}\right| = O(k^{-1}).$$
\end{theorem}

We will present a proof of this theorem here, which builds further
on results of Karabegov and Schlichenmaier \cite{KS}, which intern
again uses the Boutet de Monvel and Sj\"{o}strand expression for the
Szeg\"{o} kernel in \cite{BdMS} and then stationary phase
approximation.

We recall the part of the setting from \cite{KS} we need here.

Let $B^{(k)}_\sigma \in C^\infty(Z\times Z)$ be the Bargman
kernel. I.e. for $s\in C^\infty(M,\L^k)$, we have that
 \[\pi_\sigma^{(k)}(\psi_s)(\alpha) = \psi_{\pi_\sigma^{(k)}s}(\alpha)
 = \int_{Z} B^{(k)}_\sigma(\alpha,\beta) \psi_s(\beta)
 \tOmega(\beta).\]
In fact we get the relation that
\[B^{(k)}_\sigma(\alpha,\beta) = \langle e^{(k)}_{\beta,\sigma}, e^{(k)}_{\alpha,\sigma}\rangle
= \psi_{e^{(k)}_{\beta,\sigma}}(\alpha).\]

Since the Bargman Kernel decays faster than any power of $k$ of
the diagonal (see e.g. \cite{KS}), it follows immediately that
\begin{equation}
|\langle e^{(k)}_{\alpha_1,\sigma_0},
e^{(k)}_{\alpha_2,\sigma_1}\rangle | = O(k^{-N})\label{csip}
\end{equation}
for all $N$ and all $\alpha_i \in \L^*_{x_i}$, $x_1\neq x_2$,
$\sigma_i\in \T$, $i=1,2$. Furthermore, for all $\alpha \in Z$ and
$\sigma \in \T$ we have that
\[\left| k^{-n}|e^{(k)}_{\alpha,\sigma}| -1 \right| = O(k^{-1}).\]
In fact Zelditch \cite{Z} provides a full asymptotic expansion of
$|e^{(k)}_{\alpha_1,\sigma}|$ to all orders in $k$, but we will
not need it here.

Let us now recall the expression from \cite{KS} for the asymptotic
expansion of $B^k_\sigma$ near the diagonal.

Choose a $\sigma \in \T$ and $x_0\in M$. Let $U$ be a sufficiently
small neighbourhood of $x_0$. Let $s$ be a holomorphic frame of
$\L^*$ over $U$. Then $\alpha(x) = \frac{s(x)}{|s(x)|}$, $x\in U $
is a smooth section of $Z$ over $U$. Furthermore $\Phi(x) = log
|s(x)|$ is a potential for $\omega$, i.e. $\omega =
-i\partial\bar\partial \Phi$. By adjusting $s$, we can arrange
that $x_0$ is a stationary point for $\Phi$. Let $\tPhi\in
C^\infty(U\times \bar U)$ be an almost analytic extension of
$\Phi$ from the diagonal, in the sense of H\"{o}rmander (see \cite{KS}
for the definition of $\tPhi$). We will and may assume that
$\tPhi(y,x) = \overline{\tPhi(y,x)}$. Let $\chi\in
C^\infty(U\times \bar U)$ be given by
\[\chi(x,y) = \tPhi(x,y) - \frac12(\Phi(x) + \Phi(y))\]
Note that $\chi(x,x) = 0 $ and by Lemma 5.5 in \cite{KS} we can
assume that $\re(\chi(x,y)) <0$ for all $x\neq y$, $x,y \in U$.
Furthermore the function $y \mapsto \chi(x_0,y)$ has a
non-degenerate critical point at $y=x_0$.

\begin{theorem}[Karabegov and Schlichenmaier]\label{TKS}
There exist a unique function $b\in C^\infty(U\times U)$ such that
for any compact subset $K \subset U \times U$, there exist a
constant $C$ such that
\[\sup_{(x,y)\in K}\left| k^{-n} B^{(k)}_\sigma(\alpha(x),\alpha(y))
- e^{k\chi(x,y)}b(x,y)\right| \leq \frac Ck.\]

\end{theorem}

We remark that a full asymptotic expansion of $B^{(k)}_\sigma$ is
given in Theorem 5.6. of \cite{KS}. We also remark that $\alpha$,
$\chi$ and $b$ of course depends on $\sigma$.

\proof[Proof of Theorem \ref{TOCSA}.]

First we compute
 \begin{eqnarray*}
 \lefteqn{\langle T^{(k)}_{f,(\sigma_0,\sigma_1)} \frac{e^{(k)}_{\alpha,\sigma_0}
}{|e^{(k)}_{\alpha,\sigma_0}|}, f(x)
\frac{e^{(k)}_{\alpha,\sigma_1}
}{|e^{(k)}_{\alpha,\sigma_1}|}\rangle}\\
 & = & \frac{\langle f e^{(k)}_{\alpha,\sigma_0}, f(x) e^{(k)}_{\alpha,\sigma_1}\rangle}
 {|e^{(k)}_{\alpha,\sigma_0}||e^{(k)}_{\alpha,\sigma_1}|} \\
 & = & \frac{1}
 {|e^{(k)}_{\alpha,\sigma_0}||e^{(k)}_{\alpha,\sigma_1}|}
 \int_Z \tau^* f(\beta) \psi_{e^{(k)}_{\alpha,\sigma_0}}(\beta)
 \overline{f(x) \psi_{e^{(k)}_{\alpha,\sigma_1}}(\beta)} \tOmega(\beta)\\
 & = & \frac{1}
 {|e^{(k)}_{\alpha,\sigma_0}||e^{(k)}_{\alpha,\sigma_1}|}
 \int_Z B^{(k)}_{\sigma_0}(\beta,\alpha) B^{(k)}_{\sigma_1}(\alpha, \beta)
  \tau^* f(\beta) \overline{f(x)} \tOmega(\beta)
\end{eqnarray*}
Choose a neighbourhood $U$ as discussed above around $x$. We can
disregard the part of the integral, which is outside a compact
neighbourhood of $x$ inside $U$, because of the estimate
(\ref{csip}). The remaining integral we evaluate using the
stationary phase method (see e.g. \cite{MS}). In fact the
situation is very closed to the one consider in \cite{KS} in
formula (5.22) and the following paragraphs. We use Theorem
\ref{TKS} to change the integral to an oscillatory integral, where
the phase is expressed in terms of sums of $\chi's$ (one for each
of the $B^{(k)}$'s). The above discussed properties of the $\chi$
implies the need properties of the phase in order to to apply the
stationary phase method, which then yields
 \[\left|\langle T^{(k)}_{f,(\sigma_0,\sigma_1)} \frac{e^{(k)}_{\alpha,\sigma_0}
}{|e^{(k)}_{\alpha,\sigma_0}|}, f(x)
\frac{e^{(k)}_{\alpha,\sigma_1}
}{|e^{(k)}_{\alpha,\sigma_1}|}\rangle - |f(x)|^2 \right| = O(\frac
1k).\] Then we compute
 \begin{eqnarray*}
 \lefteqn{\langle T^{(k)}_{f,(\sigma_0,\sigma_1)} \frac{e^{(k)}_{\alpha,\sigma_0}
}{|e^{(k)}_{\alpha,\sigma_0}|},T^{(k)}_{f,(\sigma_0,\sigma_1)}
\frac{e^{(k)}_{\alpha,\sigma_0}
}{|e^{(k)}_{\alpha,\sigma_0}|}\rangle}\\
 & = &\frac{\langle f e^{(k)}_{\alpha,\sigma_0}, \pi_{\sigma_1}^{(k)}
 f e^{(k)}_{\alpha,\sigma_0}\rangle}
 {|e^{(k)}_{\alpha,\sigma_0}||e^{(k)}_{\alpha,\sigma_0}|} \\
 & = & \frac{1}
 {|e^{(k)}_{\alpha,\sigma_0}||e^{(k)}_{\alpha,\sigma_0}|}
 \int_Z \tau^* f(\beta) \psi_{e^{(k)}_{\alpha,\sigma_0}}(\beta)
 \overline{\psi_{\pi_{\sigma_1}^{(k)}
 f e^{(k)}_{\alpha,\sigma_0}}(\beta)} \tOmega(\beta)\\
 & = &\frac{1}
 {|e^{(k)}_{\alpha,\sigma_0}||e^{(k)}_{\alpha,\sigma_0}|}
 \int_{Z\times Z} \tau^* f(\beta) \psi_{e^{(k)}_{\alpha,\sigma_0}}(\beta)
   \overline{B^{(k)}_{\sigma_1}(\beta, \gamma)\tau^*f (\gamma)
   \psi_{e^{(k)}_{\alpha,\sigma_0}}(\gamma)} \tOmega(\beta)\tOmega(\gamma)\\
 & = &\frac{1}
 {|e^{(k)}_{\alpha,\sigma_0}||e^{(k)}_{\alpha,\sigma_0}|}
 \int_{Z\times Z} \tau^* f(\beta) \overline{\tau^*f (\gamma)}
 B^{(k)}_{\sigma_0}(\beta,\alpha)
   B^{(k)}_{\sigma_1}(\gamma, \beta)
   B^{(k)}_{\sigma_0}(\alpha,\gamma) \tOmega(\beta)\tOmega(\gamma)
\end{eqnarray*}
Which we treat by the same method as above (again in parallel to
formula (5.22) and the following paragraphs of \cite{KS}) to get
that
\[\left|\langle T^{(k)}_{f,(\sigma_0,\sigma_1)} \frac{e^{(k)}_{\alpha,\sigma_0}
}{|e^{(k)}_{\alpha,\sigma_0}|},T^{(k)}_{f,(\sigma_0,\sigma_1)}
\frac{e^{(k)}_{\alpha,\sigma_0}
}{|e^{(k)}_{\alpha,\sigma_0}|}\rangle - |f(x)|^2\right| = O(\frac
1k).\] But then we can simply compute
 \begin{eqnarray*}
 \lefteqn{ | T^{(k)}_{f,(\sigma_0,\sigma_1)} \frac{e^{(k)}_{\alpha,\sigma_0}
}{|e^{(k)}_{\alpha,\sigma_0}|}-
f(x)\frac{e^{(k)}_{\alpha,\sigma_1}
}{|e^{(k)}_{\alpha,\sigma_1}|}|}\\
 & = & | T^{(k)}_{f,(\sigma_0,\sigma_1)} \frac{e^{(k)}_{\alpha,\sigma_0}
}{|e^{(k)}_{\alpha,\sigma_0}|} |^2 + |f(x)|^2\\
 & &\quad - 2 \re( \langle T^{(k)}_{f,(\sigma_0,\sigma_1)} \frac{e^{(k)}_{\alpha,\sigma_0}
}{|e^{(k)}_{\alpha,\sigma_0}|},f(x)\frac{e^{(k)}_{\alpha,\sigma_1}
}{|e^{(k)}_{\alpha,\sigma_1}|}\rangle \\
 & & = O(\frac 1k).
\end{eqnarray*}

\eproof

\begin{theorem}\label{BMSG}
For all $f\in \C^\infty(M)$  and any pair of points $(\sigma_0,\sigma_1)$ in $\T$
we have that
$$ \lim_{k\ra \infty}\| T^{(k)}_{f,(\sigma_0,\sigma_1)} \| = \sup_{x\in M} |f(x)|.$$
\end{theorem}

We remark that this theorem is a generalization of Theorem \ref{BMS1}
by Bordemann, Meinrenken and Schlichenmaier.

\proof It is clear that $\| T^{(k)}_{f,(\sigma_0,\sigma_1)} \|
\leq \sup_{x\in M} |f(x)|$. Choose an $x_0\in M$ such that $
|f(x_0)|=\sup_{x\in M} |f(x)|$. Pick an $\alpha_0\in \L_{x_0}^*$
and use Theorem \ref{TOCSA} to conclude that $\|
T^{(k)}_{f,(\sigma_0,\sigma_1)} \| \geq \sup_{x\in M} |f(x)|$.
\eproof

Suppose we have a smooth section $X\in C^\infty(M, T_\sigma)$ of
the holomorphic tangent bundle of $M_\sigma$. We then claim that
the operator $\pi_\sigma \nabla_X$ is a zero-order Toeplitz
operator. Suppose $s_1\in C^\infty(M,\L^k)$ and $s_2 \in
H^0(M_\sigma,\L^k)$, then we have that
\[X(s_1,s_2) = (\nabla_X s_1, s_2).\]
Now, calculating the Lie derivative along $X$ of
$(s_1,s_2)\omega^m$ and using
 the above,
one obtains after integration that
\[\langle \nabla_X s_1, s_2 \rangle = - \langle \Lambda d(i_X\omega) s_1, s_2 \rangle,\]
where $\Lambda$ denotes contraction with $\omega$. Thus
\begin{equation}\pi_\sigma \nabla_X = T_{f_X, \sigma}^{(k)},\label{1to0order}
\end{equation}
as operators from $C^\infty(M,\L^k)$ to $H^0(M_\sigma,\L^k)$,
where $f_X = -\Lambda d(i_X\omega)$. Iterating this, we find for
all $X_1,X_2 \in C^\infty(T_\sigma)$ that
\begin{equation}\pi_\sigma \nabla_{X_1}\nabla_{X_2} = T^{(k)}_{f_{X_2}f_{X_1} - X_2(f_{X_1}),\sigma}\label{2to0order}
\end{equation}
again as operators from $C^\infty(M,\L^k)$ to
$H^0(M_\sigma,\L^k)$.

For $X\in C^\infty(M, TM\otimes \bC)$ and for $s_1, s_2 \in
C^\infty(M,\L^k)$, we have that
\[\bar X (s_1,s_2) = (\nabla_{\bar X} s_1, s_2) + ( s_1, \nabla_X s_2).\]
Computing the Lie derivative along $\bar X$ of $(s_1,s_2)\omega^m$
and integrating, we get that
\[\langle \nabla_{\bar X} s_1, s_2 \rangle + \langle (\nabla_{X})^* s_1, s_2 \rangle
 = - \langle \Lambda d(i_{\bar X}\omega) s_1, s_2 \rangle.\]
Hence we see that
\begin{equation}\label{dualnab}
(\nabla_{X})^* = - \left( \nabla_{\bar X} - f_{\bar X} \right)
\end{equation}
as operators on $C^\infty(M,\L^k)$. In particular if $X\in
C^\infty(M, T_\sigma)$, we see that
\begin{equation}\label{1to0order*}
\pi_\sigma (\nabla_{X})^* \pi_\sigma = - T_{f_{\bar
X},\sigma}|_{H^0(M_\sigma,\L^k)} : H^0(M_\sigma,\L^k) \ra
H^0(M_\sigma,\L^k).
\end{equation}
For two smooth sections $X_1, X_2$ of the holomorphic tangent
bundle $T_\sigma$ and a smooth function $h\in C^\infty(M)$, we
deduce from the formula for $(\nabla_{X})^*$ that
\begin{eqnarray}\label{2to0order*}
\pi_\sigma (\nabla_{{X}_1})^*(\nabla_{{X}_2})^* h\pi_\sigma & = &
\pi_\sigma \bar X_1 \bar X_2(h)\pi_\sigma +\\ & & \quad \pi_\sigma
f_{{\bar X}_1} \bar X_2(h)\pi_\sigma +\pi_\sigma f_{{\bar
X}_2}\bar X_1(h) \pi_\sigma +
 \nonumber \\
& & \quad \pi_\sigma\bar X_1(f_{{\bar X}_2}) h \pi_\sigma +
\pi_\sigma f_{{\bar X}_1} f_{{\bar X}_2}h \pi_\sigma\nonumber
\end{eqnarray}
as operators on $H^0(M_\sigma,\L^k)$.

Suppose $X\in \C^\infty(M,TM\otimes \bC)$. Since we have that
$T_{\sigma_0}\cap \bar T_{\sigma_1} = \{0\}$, we get a
decomposition $$ X = X' + X''$$ where $X' \in
\C^\infty(M,T_{\sigma_0})$ and $X'' \in \C^\infty(M,\bar
T_{\sigma_1})$. We then have by formula (\ref{dualnab})
$$(\nabla_X)^* = - (\nabla_{\bar X'} + \nabla_{\bar X''} - f_{\bar
X}).$$ From which we conclude for any $h\in \C^\infty(M)$ that
\[
\begin{split}
\pi_{\sigma_1} (\nabla_X)^* h \pi_{\sigma_0}& = - \pi_{\sigma_1}
(h(\nabla_{\bar X''} - f_{\bar X}) + \bar X(h))\pi_{\sigma_0}\\ &
= \pi_{\sigma_1} (f_{h\bar X''}-f_{\bar X}h + \bar
X(h))\pi_{\sigma_0}.
\end{split}
\]
Suppose we now have $X_1, X_2\in \C^\infty(M,TM\otimes \bC)$. Then
we compute
\[
\begin{split}
\pi_{\sigma_1} (\nabla_{X_1})^* (\nabla_{X_2})^* h \pi_{\sigma_0}
& = \pi_{\sigma_1} (\nabla_{\bar X_1} - f_{\bar X_1})
(\nabla_{\bar X_2} - f_{\bar X_2}) h \pi_{\sigma_0} \\ & =
\pi_{\sigma_1} \nabla_{\bar X_1} \nabla_{\bar X_2} h
\pi_{\sigma_0} \\ & \quad - \pi_{\sigma_1} f_{\bar X_1}
\nabla_{\bar X_2} h \pi_{\sigma_0} -
 \pi_{\sigma_1}   \nabla_{\bar X_1} f_{\bar X_2}h \pi_{\sigma_0}\\
& \quad + \pi_{\sigma_1} f_{\bar X_1} f_{\bar X_2}h
\pi_{\sigma_0}\\ & = - \pi_{\sigma_1} h \nabla_{\bar X_1}
\nabla_{\bar X''_2}\pi_{\sigma_0}\\ & \quad - \pi_{\sigma_1} h
f_{\bar X_1} \nabla_{\bar X''_2} \pi_{\sigma_0} -
 \pi_{\sigma_1}   h f_{\bar X_2} \nabla_{\bar X''_1}  \pi_{\sigma_0}\\
 & \quad + \pi_{\sigma_1}  \bar X_1(h) \nabla_{\bar X''_2}  \pi_{\sigma_0}  +
 \pi_{\sigma_1}   \bar X_2(h) \nabla_{\bar X''_1}  \pi_{\sigma_0}\\
  & \quad + \pi_{\sigma_1}  \bar X_1\bar X_2(h) \pi_{\sigma_0}  -
 \pi_{\sigma_1}  f_{\bar X_1} \bar X_2(h)   \pi_{\sigma_0} -
 \pi_{\sigma_1}  f_{\bar X_2} \bar X_1(h)   \pi_{\sigma_0}\\
 & \quad -\pi_{\sigma_1}  \bar X_1(f_{\bar X_2} ) h   \pi_{\sigma_0}
+ \pi_{\sigma_1} f_{\bar X_1} f_{\bar X_2}h \pi_{\sigma_0}
\end{split}
\]

Now
\[
\begin{split}
\pi_{\sigma_1} h\nabla_{\bar X_1} \nabla_{\bar X''_2}
\pi_{\sigma_0} & = \pi_{\sigma_1} h(\nabla_{\bar X'_1}
\nabla_{\bar X''_2} + \nabla_{\bar X''_1} \nabla_{\bar
X''_2})\pi_{\sigma_0} \\ & = k \pi_{\sigma_1} h \omega(\bar X'_1,
\bar X''_2) \pi_{\sigma_0} \\ & \quad \pi_{\sigma_1} h
\nabla_{[\bar X'_1, \bar X''_2]} \pi_{\sigma_0} + \pi_{\sigma_1}
h\nabla_{\bar X''_1} \nabla_{\bar X''_2} \pi_{\sigma_0}
\end{split}
\]
Hence, by splitting $[\bar X'_1, \bar X''_2]
\in\C^\infty(M,TM\otimes \bC)$ with respect to the direct sum
$TM\otimes \bC = \bar T_{\sigma_0} \oplus T_{\sigma_1}$ and by
using formula (\ref{1to0order}) and (\ref{2to0order}), we see
there exists functions $H_0(X_1,X_2)_{(\sigma_0,\sigma_1)}(h),
H_1(X_1,X_2)_{(\sigma_0,\sigma_1)}(h) \in \C^\infty(M)$ such that

\begin{eqnarray}
\lefteqn{\pi_{\sigma_1} (\nabla_{X_1})^* (\nabla_{X_2})^* h
\pi_{\sigma_0} =} \label{2orderopto}\\ & k \pi_{\sigma_1}
H_0(X_1,X_2)_{(\sigma_0,\sigma_1)}(h) \pi_{\sigma_0} +
\pi_{\sigma_1} H_1(X_1,X_2)_{(\sigma_0,\sigma_1)}(h)
\pi_{\sigma_0}.\nonumber
\end{eqnarray}
In fact $$ H_0(X_1,X_2)_{(\sigma_0,\sigma_1)}(h) = h \omega(\bar
X'_1, \bar X''_2). $$ and $H_1(X_1,X_2)_{(\sigma_0,\sigma_1)}$ is
a second order differential operator as a function of $h$.

These calculations will be applied in section \ref{cscc} in the
proof of Proposition \ref{Mainnorm}.

\section{Hermitian structures on $\V_k$ and $\End(\V_k)$.}\label{Hermstruc}

In this section we consider a further geometric Hermitian structures on
$\V_k$ and recall its asymptotic flatness as proved in \cite{A2}.

We will use the following Hermitian structure on $\cH_k$
\begin{equation}
\langle s_1,s_2\rangle_F = \frac{1}{m!}\int_M (s_1,s_2) e^{- F}\omega^m, \label{asymphermform}
\end{equation}
where $s_1, s_2$ are sections of $\cH_k$ over $\T$. The associated operator norm
on sections of $\End(\V_k)$ is denoted $\|\cdot\|_F$.

We recall that $F(\s)$ is the Ricci potential on $M_\s$
for each $\s\in \T$ determined by equation (\ref{riccipot}).
From \cite{A2} we recall that

\begin{lemma}\label{equinormL2}
The Hermitian structure on $\cH_k$
$$\langle s_1,s_2\rangle_F = \frac{1}{m!}\int_M (s_1,s_2)
e^{- F}\omega^m $$
and the constant $L_2$-Hermitian structure on $\cH_k$
$$\langle s_1,s_2 \rangle = \frac{1}{m!}\int_M (s_1,s_2) \omega^m$$
are equivalent uniformly in $k$ when restricted to $\V_k$
over any compact subset $K$ of $\T$.
\end{lemma}

The constant $L_2$-Hermitian structure on $\cH_k$ is
not asymptotically flat with respect to $\Nabla$.

For any tangent vector field $V$ on $\T$ we have
\[V\langle s_1,s_2\rangle_F =
\langle \Nabla^t_{V} s_1,s_2\rangle_F + \langle
s_1,\Nabla^t_{V} s_2\rangle_F - \langle
V[F] s_1,s_2\rangle_F.\] Let now
\[E(V)(s) = V |s|^2_F -
\langle \Nabla_{V}s,s\rangle_F -
\langle s,\Nabla_{V}s\rangle_F .\]
Then we have by Proposition 2 in \cite{A2} that
\begin{proposition}\label{asympherm}
The Hermitian structure (\ref{asymphermform}) is asymptotically flat
with respect to the connections $\Nabla$, i.e. for any compact subset $K$ of
$\T$ and any vector field $V$ defined over $K$, there exists a
constant $C$ such that for all sections $s$ of $\V_k$ over $K$, we
have that
\[|E(V)(s)| \leq \frac{C}{k+n} |s|_F^2\]
over $K$.
\end{proposition}

This Proposition is used in section \ref{cscc}.

\section{Asymptotic flatness of the coherent states.}\label{cscc}

Let us first consider the asymptotics of the parallel transport of the
Hitchin connection. Let $J= [0,1]$.

\begin{theorem}\label{Asympflat}
Let $\sigma : J \ra \T$ be a curve in \Teim space from $\sigma_0$ to $\sigma_1$ and
$P_{\sigma}$ the parallel transport in the
bundle $\V_k$ with respect to $\Nabla$ along $\sigma$.
Then there exist a function $g_\sigma\in C^\infty(M)$ such that
\[\|P_{\sigma} - T^{(k)}_{g_\sigma, (\sigma_0,\sigma_1)}\| = O(k^{-1}),\]
where $\|\cdot\|$ is the operator norm with respect to the
$L_2$-norm on $H^0(M_{\s_i},\L_{\s_i}^k)$, $i=0,1$.
\end{theorem}

In order to prove this theorem we first prove

\begin{proposition}\label{Mainnorm}
Let $\sigma : J \ra \T$ be a curve in \Teim space from $\sigma_0$ to $\sigma_1$.
Then there exist a unique curve $g : J \ra C^\infty(M)$ such that
$g(0) =1 $ and
\begin{equation}
\sup_{t\in J}|\Nabla_{\sigma'(t)}(T_{g(t),\sigma(t)}s)| = \frac Ck
|s|\label{covconstm}
\end{equation}
for some constant $C$ and all $s\in H^0(M_{\sigma_0},\L^k).$
\end{proposition}

\proof

We recall from \cite{A2} the formula
\begin{equation}
\pi_{\sigma(t)} (\pi_{\sigma(t)})' = \pi_{\sigma(t)} u(\sigma'(t))^* -
\pi_{\sigma(t)} u(\sigma'(t))^*\pi_{\sigma(t)}. \label{derivepi}
\end{equation}

For any choice of $g : J \ra C^\infty(M)$ and $s\in H^0(M_{\sigma_0},\L^k)$,
we compute using this formula that
\[
\begin{split}
\pi_{\sigma(t)} \Nabla_{\sigma'(t)}(T_{g(t),\sigma(t)} s) & =
\pi_{\sigma(t)} g'(t) \pi_{\sigma_0}s\\
            & \quad \pi_{\sigma(t)}u(\sigma'(t))^* g(t) \pi_{\sigma_0}s -
            \pi_{\sigma(t)}u(\sigma'(t))^* \pi_{\sigma(t)}g(t) \pi_{\sigma_0}s \\
            & \quad - \pi_{\sigma(t)}u(\sigma'(t)) \pi_{\sigma(t)}g(t) \pi_{\sigma_0}s.
\end{split}\]

By applying Theorem \ref{BMSG} and formulae (\ref{1to0order}), (\ref{2to0order}), (\ref{1to0order*}), (\ref{2to0order*}) and
(\ref{2orderopto}), we see that there exist a unique smooth map
$A : J \ra \cD(M)$ and a constant $C$ such that
\begin{equation}\label{evT}
| \pi_{\sigma(t)} \Nabla_{\sigma'(t)}(T_{g(t),\sigma(t)} s) -
\pi_{\sigma(t)} (g'(t) - A(t)(g(t)) (s)) | \leq \frac Ck |s|.
\end{equation}
Let now $g : J \ra \C^\infty(M)$ be the unique smooth map, such
that $g(0) = 1$ and which solves
\begin{equation}\label{ODE}
g'(t) = A(t)(g(t)).
\end{equation}
We then have the required estimate.
Conversely any smooth $g: J \ra \C^\infty(M)$ with the property that the norm
estimate (\ref{covconstm}) is satisfied most solve (\ref{ODE}) by
formula (\ref{evT}) and Theorem \ref{BMSG}. Hence
the uniqueness of $g$ follows from the uniqueness of solutions to
ordinary differential equations.
\eproof

\proof[Proof of Theorem \ref{Asympflat}]
Let $\sigma: J \ra \T$  be a smooth  curve in
$\T$ from $\sigma_0$ to $\sigma_1$. Let $g : J \ra C^\infty(M)$ be
given by Proposition \ref{Mainnorm}.
Define
\[\Theta_k(t) : \V_{k, \sigma_0} \ra \V_{k, \sigma(t)}\]
by
\[ \Theta_k(t) = P_{\sigma,t} - T^{(k)}_{g(t),\sigma(t)},\]
where $P_{\sigma,t}$ is the parallel transport in $\V_k$ along
$\sigma$ from $\sigma_0$ to $\sigma(t)$.
Let $s\in H^0(M_{\sigma_0},\L^k)$ such that $|s|_F = 1$ and define $n_k : J \ra [0,\infty)$ by
\[n_k(t) = |\Theta_k(t)s|_F^2.\]
Then the functions $n_k$ are differentiable in $t$ and we compute that
\[
\begin{split}
\frac{d n_k}{d t} & = \langle
\Nabla_{\sigma'(t)}(\Theta_k(t)s),\Theta_k(t)s\rangle_F
+ \langle
\Theta_k(t),
\Nabla_{\sigma'(t)}(\Theta_k(t)s)\rangle_F + E(\Theta_k(t)s)\\
&=- \langle \Nabla_{\sigma'(t)}T^{(k)}_{g(t),\sigma(t)}s,
\Theta_k(t)s\rangle_F - \langle
\Theta_k(t)s,
\Nabla_{\sigma'(t)} T^{(k)}_{g(t),\sigma(t)}s\rangle_F + E(\Theta_k(t)s).
\end{split}
\]
Using the above, we get the following estimate
\[
\begin{split}
|\frac{d n_k}{d t}| &\leq 2 |\Nabla_{\sigma'(t)}
T^{(k)}_{g(t),\sigma(t)}s|_F |\Theta_k(t)s|_F + |E(\Theta_k(t)s)|\\
& \leq 2 |\Nabla_{\sigma'(t)}
T^{(k)}_{g(t),\sigma(t)}s|_F n_k^{1/2} + |E(\Theta_k(t)s)|.\end{split}
\]
Now we use Lemma \ref{equinormL2} and Propositions \ref{asympherm} and \ref{Mainnorm}
to obtain that there
exists a constant $C$ independent of $s$ such that
\[
|\frac{d n_k}{d t}| \leq \frac{C}{k} ( n_k^{1/2} + n_k).
\]
This estimate implies that
\[n_k(t) \leq (\exp(\frac{C t}{2k}) - 1 )^2.\]
Thus
\[|P_{\sigma}s - T_{g(1),\sigma_1}^{(k)}s|_F
= |\Theta_k(1)s|_F\leq  C_1 n_k(1)^{1/2}.\]
The Theorem then follows from this estimates combined with
Lemma \ref{equinormL2}.
\eproof

Let $x \in M$ be a point in the moduli space. Let $\alpha$ be a
point in $\L_x-\{0\}$. As in section \ref{sec4} we consider the
associated section $e_\alpha^{(k)}$ of $\V_k$ determined by
$\alpha$ by formula (\ref{defcs}).

\begin{theorem}\label{PCS}
Let $\sigma : J \ra \T$ be a curve from $\sigma_0$ to $\sigma_1$.
Let $g_\sigma \in C^\infty(M)$ be the function determined by
$\sigma$ in Theorem \ref{Asympflat}. Then
\[|P_\sigma(\frac{e_{\alpha, \sigma_0}^{(k)}}{|e_{\alpha,
\sigma_0}^{(k)}|}) - g_\sigma(x)\frac{e_{\alpha,
\sigma_0}^{(k)}}{|e_{\alpha, \sigma_0}^{(k)}|}| = O(k^{-1}).\]
\end{theorem}

This theorem follows directly from Theorem \ref{Asympflat} and
Theorem \ref{TOCSA}.

Define a section $E_x^{(k)}$ of $\End(\V_k)$ from the section
$e_\alpha^{(k)}$ of $\V_k$ as follows
\[E_x^{(k)}(s) = \langle s, e_\alpha^{(k)}\rangle e_\alpha^{(k)}. \]

\proof[Proof of Theorem \ref{ccE}]
Let $\sigma : J \ra \T$ be a curve from $\sigma_0$ to $\sigma_1$.
For $s\in H^0(M_{\sigma_1},\L^k)$, we have that
\[
\begin{split}
P^e_\sigma(E^{(k)}_{x}(\sigma_0))(s) & = P_\sigma \circ
E^{(k)}_x(\sigma_0) P_\sigma^{-1}(s)\\ & = \langle
P_\sigma^{-1}(s),\frac{e_{\alpha, \sigma_0}^{(k)}}{|e_{\alpha,
\sigma_0}^{(k)}|} \rangle P_\sigma(\frac{e_{\alpha,
\sigma_0}^{(k)}}{|e_{\alpha, \sigma_0}^{(k)}|})\\ & = \langle
s,(P_\sigma^{-1})^*(\frac{e_{\alpha, \sigma_0}^{(k)}}{|e_{\alpha,
\sigma_0}^{(k)}|})\rangle P_\sigma(\frac{e_{\alpha,
\sigma_0}^{(k)}}{|e_{\alpha, \sigma_0}^{(k)}|}). \end{split}\] But
then by theorem \ref{PCS} we get that
\[|P^e_\sigma(E^{(k)}_{x}(\sigma_0)) - \langle\cdot, \overline{g_\sigma(x)}^{-1}\frac{e_{\alpha, \sigma_1}^{(k)}}{|e_{\alpha,
\sigma_1}^{(k)}|}\rangle \otimes g_\sigma(x) \frac{e_{\alpha,
\sigma_1}^{(k)}}{|e_{\alpha, \sigma_1}^{(k)}|}| = O(k^{-1}).\]
Which implies the theorem since $$ E^{(k)}_{x}(\sigma_1) =
\langle\cdot, \overline{g_\sigma(x)}^{-1}\frac{e_{\alpha,
\sigma_1}^{(k)}}{|e_{\alpha, \sigma_1}^{(k)}|}\rangle \otimes
g_\sigma(x) \frac{e_{\alpha, \sigma_1}^{(k)}}{|e_{\alpha,
\sigma_1}^{(k)}|}.$$ \eproof

\section{The almost fixed vector}\label{PF}

Let $\Lambda \subset SU(2)$ be the finite subgroup
\[\Lambda = \langle \left( \begin{array}{cc}
  i & 0 \\
  0 & -i
\end{array} \right), \left( \begin{array}{cc}
  0 & 1 \\
  -1 & 0
\end{array} \right)\rangle.\]

Let $X$ be the image
\[X = \im ( \Hom'(\tilde \pi_1(\Sigma), \Lambda) \ra M).\]
Choosing a standard set of
generators for $\pi_1(\Sigma - \{p\})$ and mapping the $i'$th pair
to the two generators of $\Lambda$ given above and the rest to the
identity, we see that $X$ contains at least
$g$ points. Moreover, $X$ is
clearly a $\Gamma$-invariant subset of $M$.

As stated in the introduction, we then define the section $E_X^{(k)}$
of $\End(\V_k)$ by
\[E_X^{(k)} = \sum_{x\in X} E_x^{(k)}\]
and let $E_{X,0}^{(k)}$ be the traceless part of
$E_X^{(k)}$. Thus
\[E_{X,0}^{(k)} = E_{X}^{(k)} - \frac{\tr(E_{X}^{(k)})}{d_g(k)} \Id.\]
By formula (\ref{csip}) we have that $\tr(E_{X}^{(k)}(\sigma_0))$
converges to $|X|$ and therefore so does $$|E_{X,0}^{(k)}|^2 =
\langle E_{X,0}^{(k)}(\sigma_0), E_{X,0}^{(k)}(\sigma_0)\rangle =
\tr((E_{X,0}^{(k)}(\sigma_0))^2),$$ since $d_g(k)$ is a polynomial
of degree $3g-3$. Hence we see that for sufficiently large $k$,
$E_{X,0}^{(k)}\neq 0$.

Hence for large enough $k$ we have a unique unit
vector $\E_{X,0}^{(k)}\in \cH_k$, which at $\sigma_0$ agrees with
$E_{X,0}^{(k)}(\sigma_0)/[E_{X,0}^{(k)}(\sigma_0)]$.

\proof[Proof of Theorem \ref{AFV}]
Suppose we have a $\phi \in \Gamma$. We then have that
\[\phi^*(E_{X,0}^{(k)}) = E_{X,0}^{(k)}\]
which means that
\[\phi^*\circ E_{X,0}^{(k)}(\sigma_0) \circ (\phi^*)^{-1} = E_{X,0}^{(k)}(\phi(\sigma_0)).\]
Hence
\[
\begin{split}
\phi(\E_{X,0}^{(k)}) & = \frac1{[E_{X,0}^{(k)}(\sigma_0)]}
P^e_{\phi(\sigma_0),\sigma_0}( E_{X,0}^{(k)}(\phi(\sigma_0)))\\
& = \frac1{[E_{X,0}^{(k)}(\sigma_0)]}
(P^e_{\phi(\sigma_0),\sigma_0}( E_{X}^{(k)}(\phi(\sigma_0))) - \frac{\tr(E_{X}^{(k)})}{d_g(k)} \Id).
\end{split}\]
But from Theorem \ref{normeq} and Theorem \ref{ccE} we then get that there exist a
constant $\tilde C$ such that
\[[\E_{X,0}^{(k)} - \phi(\E_{X,0}^{(k)})] \leq
\frac C{[E_{X,0}^{(k)}(\sigma_0)]}
|E_{X}^{(k)}(\sigma_0) - P^e_{\phi(\sigma_0),\sigma_0}
( E_{X}^{(k)}(\phi(\sigma_0)))| \leq \frac{\tilde C}k\]
for all sufficiently large $k$.
\eproof

\end{document}